\documentclass[10pt, a4paper]{amsart}
\usepackage[charter]{mathdesign}
\usepackage{enumitem}
\usepackage{amsmath, wrapfig}
\usepackage{fancybox, graphicx, subfigure}
\usepackage[usenames, dvipsnames, pdftex]{color}
\usepackage[colorlinks=true,
        raiselinks=true,
        linkcolor=MidnightBlue,
        citecolor=Mahogany,
        urlcolor=ForestGreen,
        pdfauthor=,
        pdftitle={},
        pdfkeywords={},
        pdfsubject={Technical Report},
        plainpages=false]{hyperref}

\allowdisplaybreaks

\numberwithin{equation}{section}

\newtheorem{theorem}{Theorem}
\newtheorem{corollary}[theorem]{Corollary}
\newtheorem{lemma}[theorem]{Lemma}
\newtheorem{proposition}[theorem]{Proposition}

\newtheorem{definition}[theorem]{Definition}

\renewcommand{\le}{\leqslant}
\renewcommand{\ge}{\geqslant}
\renewcommand{\leq}{\leqslant}
\renewcommand{\geq}{\geqslant}

\newcommand{\Char}[1]{\boldsymbol{1}_{#1}}
\newcommand{\indic}[1]{\boldsymbol{1}_{#1}}

\newcommand{\PP}{\ensuremath{\mathbb P}}
\newcommand{\EE}{\ensuremath{\mathbb E}}
\newcommand{\sigalg}{\ensuremath{\mathfrak F}}
\newcommand{\Borelsigalg}[1]{\ensuremath{\mathfrak B}(#1)}
\newcommand{\R}{\mathbb R}
\newcommand{\N}{\mathbb N}
\newcommand{\Nz}{{\N_0}}

\newcommand{\mrm}{\mathrm}

\newcommand{\mn}{\wedge}

\newcommand{\setmin}{\setminus}
\newcommand{\Lop}{L}

\title{On the connections between PCTL and Dynamic Programming}

\author[F.~Ramponi]{Federico Ramponi}
\author[D.~Chatterjee]{Debasish Chatterjee}
\author[S.~Summers]{Sean Summers}
\author[J.~Lygeros]{John Lygeros}
\email{\{ramponif,chatterjee,ssummers,lygeros\}@control.ee.ethz.ch}
\thanks{This research was partially supported by the Swiss National Science Foundation, grant 200021-122072.}

\address{Automatic Control Laboratory, ETL I28, ETH Z\"urich, Physikstrasse 3, 8092 Z\"urich, Switzerland}
\urladdr{\url{http://control.ee.ethz.ch}}


\begin{document}
\maketitle


\begin{abstract}
Probabilistic Computation Tree Logic
(PCTL) is a well-known modal logic which has become a standard
for expressing temporal properties of finite-state Markov chains
in the context of automated model checking.
In this paper, we give a definition of PCTL for {\em noncountable}-space
Markov chains, and we show that there is a substantial affinity
between certain of its operators and problems of Dynamic Programming.
After proving some uniqueness properties of the solutions to the latter,
we conclude the paper with two examples to show that some
recovery strategies in practical applications, which are naturally stated
as reach-avoid problems, can be actually viewed as particular cases of PCTL 
formulas. 
\end{abstract}


\section{Introduction}
\label{SECTION_INTRODUCTION}

Reachability analysis of deterministic dynamical systems constitutes a practically important and intensely researched area in control theory.  Over the years, a wide variety of tools and methods have been developed to verify the dynamic properties of these systems, for examples see \cite{Lygeros1, LygCont, Abate2007, AminSHS, Aubin1, Mitch}.  In particular, in \cite{Lygeros1, Mitch} the reachability problems considered are solved via dynamic programming (DP).  As a result, a large number of exact and approximate methods for solving the central \emph{Bellman equation} in DP \cite{ref:bertsekasDP1, ref:bertsekasDP2, ref:bertsekasNDP, ref:powellADP} can be exploited for the solution of verification problems of deterministic dynamical systems.  

Recently, reachability analysis of stochastic Markovian processes has gained significant interest, and mechanisms for the verification of safety and performance properties by means of a control policy have been explored.  An example of such a problem is to find the probability, starting from a certain state $x$, of reaching a ``safe'' set within a certain number of time-steps, where the state $x$ could be labelled ``almost safe'' if such probability is greater than, say, $1-\varepsilon$.  A related problem, which has been studied recently by some of the authors, is that of {\em maximizing the probability of reaching a ``safe'' set, while avoiding a ``bad'' set}~\cite{ref:seanpreachavoid, ref:preachavoidarXiv}. This problem arose as a remedy for the impossibility of imposing hard state constraints in stochastic model predictive control.  In general, if one considers an infinite trajectory of a stochastic system, every compact state-constraint set is going to be violated almost surely at some time.  Thus, a good course of action when this happens is to devise a recovery strategy to drive the controlled system from the ``unsafe'' states back to the set of ``safe'' states.

If a control variable is unavailable or a control policy has been predetermined, the verification of the stochastic system reduces to calculating the likelihood of the occurance of certain events.  In this manner, the above problem is directly related to stochastic model checking of {\em finite}-state Markov models in that the analysis involves both reachability and likelihood computations.  Therefore, it is reasonable to consider an extension of Probabilistic Computation Tree Logic (PCTL), a modal logic developed for {\em finite}-state Markov chains, which forms the foundation for the automated verification tools for {\em finite}-state Markov models, to general state-space Markov chains.  

Algorithms for stochastic model checking {\em finite}-state Markov models come from standard deterministic model checking, linear algebra, and the analysis of Markov chains.  Finite state model checkers include the software tools PRISM \cite{PRISM}, SMART \cite{SMART}, $E\vDash MC^2$ \cite{EMC2}, and MRMC \cite{MRMC}, and have been used to solve various problems over the last few years.  In the area of systems biology, probabilistic model checking has been used in the analysis of biological pathways \cite{BIO1,BIO2} and signalling events \cite{BIO3}.  Additional examples of the use of stochastic model checking include the probabilistic verification of security protocols \cite{SEC}, dynamic power management \cite{DPM}, and residual risks in safety-critical systems \cite{RISK}.

In this paper, we consider the verification of general state-space Markov chains through an extension of the standard grammar and semantics of PCTL to non-countable-state Markov chains (the reader can find a similar extension in \cite{HUTH_ONFINITESTATE}).
As with the finite case, the evaluation of a PCTL formula can be recursively reduced to the truth of atomic propositions by employing computations dictated by the PCTL semantics. In this process of reduction, certain rules of the semantics simply stipulate unions or intersections of sets, while others involve the computation of integrals. It is in the computation of the integrals where the bulk of the algorithmic methodology is contained.  We show that the ``bounded until'' operator, which considers the property of hitting a ``safe'' from an ``unsafe'' set over a finite time horizon, can be evaluated through a dynamic recursion. Additionally, we prove that the ``unbounded until'' operator, which considers the property of hitting a ``safe'' set from an ``unsafe'' set at some point in time, can be evaluated via a DP-like \emph{Bellman equation}.  Further, we emphasize that, while in the numerical examples provided we grid the state space in order to solve the integral equations, any method in the literature for the numerical computation of a DP can be exploited for this problem.

{\em Outline of the work}:  In section \ref{SECTION_PCTL} we review the standard grammar and semantics of PCTL for finite-state Markov chains.  In section \ref{SECTION_PCTL_GENERAL} we extend the grammar and semantics of PCTL to general state-space Markov chains.  The uniqueness of a certain function associated with the ``unbounded until'' property is considered in section \ref{SECTION_STOPPING_TIMES}. Finally, section \ref{SECTION_EXAMPLES} concludes the paper with some applications and numerical examples.

\section{Probabilistic Computation Tree Logic}
\label{SECTION_PCTL}

\newcommand{\MCX}{X} 
\newcommand{\MCx}{x} 
\newcommand{\MCinitial}{x} 
\newcommand{\MCP}{Q} 
\newcommand{\MCA}{{\mathcal A}} 

\newcommand{\PCTLsetA}{\Phi}
\newcommand{\PCTLsetB}{\Psi}
\newcommand{\PCTLsetAMB}{{\PCTLsetA\backslash\PCTLsetB}}
\newcommand{\PathFormula}{\phi}
\newcommand{\Sat}{{\rm Sat}}
\newcommand{\Trajectory}{\omega}
\newcommand{\TrajectorySpace}{\Omega}

\newcommand{\PCTLP}[2]{P_{#1}\left[#2\right]}
\newcommand{\PCTLPP}[1]{\PCTLP{\sim p}{#1}}
\newcommand{\PCTLX}{\mathcal{X}}
\newcommand{\Until}{{\mathcal U}}
\newcommand{\PCTLBuntil}[3]{#1\ \Until^{\leq#2}\ #3}
\newcommand{\PCTLuntil}[2]{#1\ \Until\ #2}
\newcommand{\Eventually}{\lozenge}
\newcommand{\Always}{\square}
\newcommand{\PCTLBeventually}[2]{\Eventually^{\leq#1}\ #2}
\newcommand{\PCTLeventually}[1]{\Eventually\ #1}
\newcommand{\True}{\mathsf{T}}
\newcommand{\False}{\mathsf{F}}

\newcommand{\ifonlyif}{{\Leftrightarrow}}
\newcommand{\Product}{\prod} 

In this section we quickly review the definition and semantics of 
PCTL for finite-state Markov chains. The reader is referred to
the original paper \cite{HANSSON_JONSSON_ALOGIC} or to the excellent survey
\cite{KWIATOWSKA_STOCHASTIC_MODEL_CHECKING} for a detailed exposition.
\\

\subsection{Labelled Markov chains}

\begin{definition}
A homogeneous, discrete-time, finite-state Markov Chain is a triple $(\MCX, \MCinitial, \MCP)$,
where:
\begin{itemize}[label=$\circ$, leftmargin=*, itemsep=0pt, partopsep=0pt, topsep=0pt]
\item $\MCX$ is a finite set of states;
\item $\MCinitial$ is the initial state; 
\item $\MCP$ is a {\em transition probability matrix}, which assigns to each pair of states $(\MCx_i, \MCx_j)$ the probability $\MCP_{\MCx_i, \MCx_j}$ of going from the state $\MCx_i$ to the state $\MCx_j$ at a given time.
\end{itemize}
\end{definition}

Consider the sample space $\Omega := \Product_{i=0}^\infty \MCX$, containing
the possible trajectories $\Trajectory = (\MCx_0, \MCx_1, ..., \MCx_t, ...)$ of the chain,
and the product $\sigma$-algebra $\sigalg$ on $\Omega$.
For a given trajectory $\Trajectory = (\MCx_0, \MCx_1, ..., \MCx_t, ...)$,
let $\Trajectory(t) := \MCx_t$.
It can be shown~\cite[pp.~90-91]{ref:borkarProbTheo} that there exists a unique probability measure on $(\Omega, \sigalg)$,
denoted by $\PP_x(\cdot)$, such that $\PP_\MCx(X_0 = \MCinitial ) = 1$
and $\PP_x( X_{t+1} = \MCx_{t+1} \mid X_t = \MCx_t, X_{t-1} = \MCx_{t-1}, ..., X_0 = \MCinitial) = \MCP(\MCx_t, \MCx_{t+1})$.

\begin{definition}
\label{DEFINITION_LMC}
Let $\MCA$ be a finite set of {\em atomic propositions}.
A {\em labelled} Markov Chain is a quadruple $(\MCX, \MCinitial, \MCP, L)$,
where:
\begin{itemize}[label=$\circ$, leftmargin=*, itemsep=0pt, topsep=0pt, partopsep=0pt]
\item $(\MCX, \MCinitial, \MCP)$ is a finite-state Markov chain;
\item $L: \MCX \rightarrow 2^\MCA$ is a set-valued function that assigns
to each state $\MCx \in \MCX$ the set $L(\MCx)\subset\MCA$ of all those
atomic propositions that are true in the state. 
\end{itemize}
\end{definition}

\subsection{Grammar and semantics of PCTL}

The grammar of PCTL is as follows:
\begin{itemize}[label=$\circ$, leftmargin=*, itemsep=0pt, topsep=0pt, partopsep=0pt]
\item $\True$ is a formula (meaning ``true'').
\item Each atomic proposition in $A\in \MCA$ is a formula.
\item If $\PCTLsetA$ and $\PCTLsetB$ are formulas, then
$\lnot \PCTLsetA$ and
$\PCTLsetA \land \PCTLsetB$
are formulas.
\item If $\PathFormula$ is a ``{\em path formula}'' (see below) and $p\in[0,1]$,
then $\PCTLPP{\PathFormula}$ is a ({\em state}) formula.
Here and throughout the rest of the paper, $\sim$
is just shortand for one of the relations $<$, $\leq$, $>$, or $\geq$.
For example, $\PCTLP{\geq 0.9}{\PathFormula}$
is one such formula, where ``$\sim$''\ $\equiv$\ ``$\geq$'' and $p=0.9$.
\end{itemize}
The above grammar defines {\em state} formulas, 
that is, formulas whose truth can be decided for each
state $\MCx\in\MCX$.
The meaning of the formulas in the first three points is the 
usual one in the standard logic of propositions.
The other standard formulas and operators can be obtained by means of
combinations of the above ones. For example, 
$\False$ (``false'') can be defined as $\lnot\True$,
$\PCTLsetA \lor \PCTLsetB$ (``inclusive or'') as $\lnot (\lnot\PCTLsetA \land \lnot\PCTLsetB)$,
and $\PCTLsetA \rightarrow \PCTLsetB$ (formal implication) as $\lnot \PCTLsetA \lor \PCTLsetB)$.

The last kind of formula is what makes PCTL a {\em modal} logic,
since it allows to express the fact that, with {\em probability}
contained in some range, something will happen in {\em time}.
It relies on the definition of {\em path} formulas,
that is, formulas whose truth is decided for 
paths $\Trajectory \in \TrajectorySpace$.
A formula like $\PCTLP{\leq 0.9}{\PathFormula}$ means, intuitively, that the probability of taking
a path that satisfies $\PathFormula$ is at least $0.9$.
If $\PCTLsetA$ and $\PCTLsetB$ are {\em state} formulas,
we define the following to be {\em path} formulas: 
\begin{itemize}[label=$\circ$, leftmargin=*, itemsep=0pt, topsep=0pt, partopsep=0pt]
\item $\PCTLX\PCTLsetA$\ \  (``next'');
\item $\PCTLBuntil{\PCTLsetA}{k}{\PCTLsetB}$\ \ (``bounded until'');
\item $\PCTLuntil{\PCTLsetA}{\PCTLsetB}$\ \ (``unbounded until'').
\end{itemize}
Intuitively, 
$\PCTLX\PCTLsetA$
means that next state will satisfy $\PCTLsetA$;
$\PCTLBuntil{\PCTLsetA}{k}{\PCTLsetB}$
means that at some time $i$, within $k$ steps, $\PCTLsetB$ will become true,
and until that time $\PCTLsetA$ will remain true; and
$\PCTLuntil{\PCTLsetA}{\PCTLsetB}$
means that at some arbitrarily large time $i$, $\PCTLsetB$ will become true,
$\PCTLsetA$ being true until then. (See the semantics below for a precise definition.)

For example, the statement 
$\MCinitial \vDash \PCTLP{\geq 0.9}{\PCTLBuntil{\PCTLsetA}{10}{\PCTLsetB}}$
means:
With probability at least $0.9$, starting from the state $\MCinitial$, 
within $10$ steps $\PCTLsetB$ will become true, and until then $\PCTLsetA$ will remain true.
In a sense, the formula 
$\PCTLP{\geq 0.9}{\PCTLBuntil{\PCTLsetA}{10}{\PCTLsetB}}$
itself denotes the set of all states $\xi$ such that,
starting from $\xi$, with probability at least $0.9$, etc.
The above statement is of course equivalent to $x$ being a member of such a set.

The two ``until'' operators allow us to define other operators which are
standard in any temporal logic.
For example, given a state formula $\PCTLsetA$,
the path formula
$\PCTLBeventually{k}{\PCTLsetA}$, which means that {\em eventually,
within $k$ steps, $\PCTLsetA$ will happen}, can be defined as
$\PCTLBuntil{\True}{k}{\PCTLsetA}$, 
and the path formula
$\PCTLeventually{\PCTLsetA}$, which means that {\em eventually, at some time,
$\PCTLsetA$ will happen}, can be defined as $\PCTLuntil{\True}{\PCTLsetA}$.
Formulas containing the standard ``always'' operator $\Always$ can also be defined,
although not in the straightforward way one may expect at first sight:
$\Always\PCTLsetA := \lnot\Eventually\lnot\PCTLsetA$ is not a correct
definition, since PCTL does not allow for the negation of
path formulas. See \cite{KWIATOWSKA_STOCHASTIC_MODEL_CHECKING} for details.

Let $A$ denote an atomic proposition, $\PCTLsetA$ and $\PCTLsetB$
denote two {\em state} formulas, and $\PathFormula$ denote a {\em path} formula.
The semantics of PCTL is defined as follows:
\begin{displaymath}
\begin{split}
\MCx \vDash \True  &\quad {\rm for\ all} \quad  \MCx \in \MCX\\
\MCx \vDash A      &\quad \ifonlyif \quad  A \in L(\MCx)\\
\MCx \vDash \lnot \PCTLsetA                  &\quad \ifonlyif \quad  \MCx \nvDash \PCTLsetA \\
\MCx \vDash \PCTLsetA \land \PCTLsetB        &\quad \ifonlyif \quad  \MCx \vDash \PCTLsetA  \ \ {\rm and}\ \ \MCx \vDash \PCTLsetB\\
\MCx \vDash \PCTLPP{\PathFormula}   &\quad \ifonlyif \quad  \PP_\MCx\left( \{\PathFormula\} \right) \sim p \\
\end{split}
\end{displaymath}

With loose notation, $\{\PathFormula\}$ stands for the set of all the
paths $\Trajectory$ that satisfy a given path formula $\PathFormula$.
Here is the related semantics:
\begin{equation}
\label{SEMANTICS_PATH_FORMULAS}
\begin{split}
\Trajectory \vDash \PCTLX\PCTLsetA                      & \ifonlyif  \Trajectory(1) \vDash \PCTLsetA    \\
\Trajectory \vDash \PCTLBuntil{\PCTLsetA}{k}{\PCTLsetB} & \ifonlyif  
  \exists i \leq k: \Trajectory(i) \vDash \PCTLsetB \ \ {\rm and}\ \  \forall j<i, \  \Trajectory(j) \vDash \PCTLsetA \\
\Trajectory \vDash \PCTLuntil{\PCTLsetA}{\PCTLsetB}     & \ifonlyif 
  \exists i \in \Nz: \Trajectory(i) \vDash \PCTLsetB \ \ {\rm and}\ \  \forall j<i, \  \Trajectory(j) \vDash \PCTLsetA \\
\end{split}
\end{equation}
Due to the latter definitions, if $\PathFormula$ is a path formula
then $\{\PathFormula\}\equiv\{\Trajectory: \Trajectory\vDash\PathFormula\}$
is always an {\em event}, that is, it always belongs to $\sigalg$.

The great relevance of PCTL for finite Markov chains lies, above all,
in the fact that the validity of arbitrarily complex formulas at
a given state can be decided exactly and in finite time. In particular,
dealing with the common operators $\lnot$, $\land$, $\lor$ etc.\ requires
just the parsing of a tree of sub-formulas; a ``bounded until'' formula
can be decided recursively; and an ``unbounded until'' formula requires
the solution of a system of linear equations.
For these matters the reader is referred to 
\cite{HANSSON_JONSSON_ALOGIC} and \cite{KWIATOWSKA_STOCHASTIC_MODEL_CHECKING}.
We shall not delve into details here, 
because the relatively easy methods available for finite Markov chains
cannot be easily extended to the case of noncountable-space Markov processes,
with respect to which the decision of PCTL formulas will be a matter of computing integrals
recursively, or solving integral equations.

\section{PCTL for general Markov processes}
\label{SECTION_PCTL_GENERAL}

\newcommand{\StateSpaceX}{X}
\newcommand{\StochasticKernelsSet}{{\mathcal P}}
\newcommand{\Atomic}{{\mathcal A}}
\newcommand{\Borel}{{\mathfrak B}}
\newcommand{\MeasurableBounded}{{{\mathcal M}_b}}
\newcommand{\SolutionsSet}{{{\mathcal M}_b^{01}(\PCTLsetA, \PCTLsetB)}}

In what follows we define PCTL grammar and semantics on a noncountable space
$X$ in terms of a stochastic kernel $Q$ and a probability measure $\PP_x$ defined
on the space of trajectories of the process.
The reader is also referred to \cite{HUTH_ONFINITESTATE} for an abstract extension of PCTL to
general Markov chains.

%
Given a nonempty Borel set $X$ (i.e., a Borel subset of a Polish space), its Borel $\sigma$-algebra is denoted by $\Borelsigalg{X}$. By convention, when referring to sets or functions, ``measurable'' means ``Borel-measurable.'' If $X$ is a nonempty Borel space, a \emph{stochastic kernel} on $X$ is a map $Q:X\times\Borelsigalg X\to[0, 1]$ such that $Q(x, \cdot)$ is a probability measure on $X$ for each fixed $x\in X$, and $Q(\cdot, B)$ is a measurable function on $X$ for each fixed $B\in\Borelsigalg X$.

Let $\StateSpaceX$ be a nonempty Borel set, and let $Q(\cdot , \cdot)$ be a stochastic kernel on $\StateSpaceX$.
For each $t=0, 1, \ldots,$ we define the space $H_t$
of \emph{admissible histories} up to time $t$ as $H_t := \Product_{i=0}^t \StateSpaceX,\; t\in \Nz$.
A generic element $h_t$ of $H_t$, called an admissible $t$-history is a vector of the form
$h_t = (x_0, x_1, \ldots, x_{t})$, with $x_j\in \StateSpaceX$ for $j=0, \ldots, t$.
Hereafter we let the $\sigma$-algebra generated by the history $h_t$ be denoted by 
$\sigalg_t$, $t\in\Nz$. Suppose the initial state $x$ is given, and let $\delta_x$ denote 
the Dirac measure at $\{x\}$. We consider the canonical sample space $\Omega := \Product_{i=0}^\infty \StateSpaceX$
and the product $\sigma$-algebra $\sigalg$ on $\Omega$. By a standard result of
Ionescu-Tulcea~\cite[Chapter~4, \S3, Theorem~5]{ref:raoProbTheo} there exists a unique
probability measure, denoted by $\PP_x(\cdot)$ on the measurable space $(\Omega, \sigalg)$
such that $\PP_x(X_0\in B) = \delta_x(B)$ and $\PP_x(X_{t+1}\in B\mid h_t) = Q(x_t, B)$ for $B\in\Borelsigalg{\StateSpaceX}$.

\subsection{Grammar and semantics}

The ``labelling'' function $L$ is introduced in
\cite{HANSSON_JONSSON_ALOGIC} and \cite{KWIATOWSKA_STOCHASTIC_MODEL_CHECKING}
as a means to specify which states satisfy which atomic propositions.
In other words, it is just a particular way to look at the {\em relation}
``$x$ satisfies $A$''.
It should be clear that an equally legitimate way to accomplish the same 
is to substitute from the beginning the ``labelling'' function
$L:\MCX\to 2^\MCA$ with a function $S:\MCA\rightarrow 2^\MCX$, that assigns to each atomic proposition
$A$ the set $S(A)$ of all those states that satisfy $A$.
The semantics can be redefined accordingly in a straightforward way:
\begin{displaymath}
\MCx \vDash A  \quad \ifonlyif \quad  \MCx \in S(A)
\end{displaymath}
But since there is no substantial difference between
saying that a state $x$ satisfies a given property,
and stating that $x$ belongs to a set, namely the set of all the states
that satisfy that property, 
it is easily seen that proceeding along this way one may drop
{\em tout-court} the distinction between formulas and sets of states satisfying them.
In the following, we shall follow this idea consistently
(mainly for ease of notation).
Thus, from now on, we shall assume that the properties expressed by formulas
are actually encoded by {\em measurable} sets $\PCTLsetA \subset \StateSpaceX$, 
we will use the letters $A, \PCTLsetA, \PCTLsetB, ...$ for
both the formulas (or atomic propositions) and the sets that encode them,
and we will use the notations $x\vDash \PCTLsetA$
and $x\in \PCTLsetA$ somewhat interchangeably.
In the same fashion, we will drop the distinction between
{\em path} formulas and {\em events} in the process's probability space.

Let us denote the family of atomic propositions with a family of
Borel measurable sets $\Atomic \subset \Borel(\StateSpaceX)$, where $\StateSpaceX\in \Atomic$.
The grammar of PCTL is defined exactly as before:
\begin{itemize}[label=$\circ$, leftmargin=*, itemsep=0pt, topsep=0pt, partopsep=0pt]
\item $\True$ is a formula (encoded by the whole space $\StateSpaceX$).
\item Each atomic proposition $A\in \MCA$ is a formula.
\item If $\PCTLsetA$ and $\PCTLsetB$ are formulas, then
$\lnot \PCTLsetA$ and
$\PCTLsetA \land \PCTLsetB$
are formulas.
\item If $\PathFormula$ is a {\em path formula} and $p\in[0,1]$,
then $\PCTLPP{\PathFormula}$ is a ({\em state}) formula.
\end{itemize}
The following are {\em path} formulas: $\PCTLX\PCTLsetA$, $\PCTLBuntil{\PCTLsetA}{k}{\PCTLsetB}$, and $\PCTLuntil{\PCTLsetA}{\PCTLsetB}$.
%

Now we define the semantics of PCTL formulas for each possible initial state
$x \in \StateSpaceX$. 
Let $A$ denote an atomic proposition and 
$\PCTLsetA$ and $\PCTLsetB$ denote formulas (measurable sets).
We define:
\begin{align*}
\MCx \vDash \True  &\quad {\rm for\ all} \quad  \MCx \in \MCX\\
x \vDash A                                & \quad \ifonlyif \quad  x \in A \\
x \vDash \lnot \PCTLsetA                  & \quad \ifonlyif \quad  x \in \PCTLsetA^C\\
x \vDash \PCTLsetA \land \PCTLsetB        & \quad \ifonlyif \quad  x \in \PCTLsetA \cap \PCTLsetB\\
x \vDash \PCTLPP{\PathFormula}            & \quad \ifonlyif \quad  \PP_x\left( \{\PathFormula\} \right) \sim p \\
\end{align*}
As in the finite state case, we can also define 
$\False := \lnot\True$, $\PCTLsetA \lor \PCTLsetB := \lnot (\lnot\PCTLsetA \land \lnot\PCTLsetB)$,
and $\PCTLsetA \rightarrow \PCTLsetB := \lnot \PCTLsetA \lor \PCTLsetB)$,
and of course we have
\begin{displaymath}
\begin{split}
x \vDash \PCTLsetA \lor  \PCTLsetB        & \quad \Leftrightarrow \quad  x \in \PCTLsetA \cup \PCTLsetB\\
x \vDash \PCTLsetA \rightarrow \PCTLsetB  & \quad \Leftrightarrow \quad  x \in \PCTLsetA^C \cup \PCTLsetB
\end{split}
\end{displaymath}
Note that all the formulas obtainable from atomic propositions
by means of the operators $\lnot, \land, \lor, \rightarrow$
are encoded by sets that belong to $\sigma(\Atomic)$.
The semantics of path formulas is defined exactly
as in equation (\ref{SEMANTICS_PATH_FORMULAS}).

\subsection{``Next''}

We will now examine the state formulas derived from the three
path formulas in greater detail. The formula arising from the ``next''
operator is trivial. Indeed,
\begin{displaymath}
\PP_x\left(\PCTLX\PCTLsetA\right)
= \PP_x\left( \{\Trajectory: \Trajectory(1) \in \PCTLsetA \} \right) 
= \PP_x\left( X_1 \in \PCTLsetA \right)
= Q(x, \PCTLsetA)
\end{displaymath}
Hence,
\begin{displaymath}
x \vDash \PCTLPP{\PCTLX\PCTLsetA} \quad \Leftrightarrow \quad  Q(x, \PCTLsetA) \sim p
\end{displaymath}
Note that
$\PCTLPP{\PCTLX\PCTLsetA}$ is a measurable set in its own right. For example,
$\PCTLP{\leq 0.5}{\PCTLX\PCTLsetA}$ is the $0.5$-sub-level set of the
measurable function $Q(\cdot, \PCTLsetA)$.
Indeed, for each $\PCTLsetA \in \Borel(\StateSpaceX) $, the set
$\{x: Q(x, \PCTLsetA) \sim p\}$ belongs to $\Borel(\StateSpaceX)$ by the
measurability of $Q(\cdot, \PCTLsetA)$.

\subsection{``Bounded until''}

Suppose that the process starts from $x_0 = x$.
On the probability space of our Markov process we define the following event:
\begin{equation}
\begin{aligned}
&\left\{ \PCTLBuntil{\PCTLsetA}{k}{\PCTLsetB} \right\}
:= \{x\in \PCTLsetB\} \cup \{x \in\PCTLsetA, x_1 \in \PCTLsetB\} \cup\\
& \quad \{x, x_1 \in \PCTLsetA, x_2 \in \PCTLsetB\} \cup ... \cup \{x, x_1, ..., x_{k-1} \in \PCTLsetA, x_k \in \PCTLsetB\}  \\ 
&= \{x\in \PCTLsetB\} \sqcup \{x \in\PCTLsetAMB, x_1 \in \PCTLsetB\} \sqcup\{x, x_1 \in \PCTLsetAMB, x_2 \in \PCTLsetB\}\\
& \qquad\sqcup ... \sqcup \{x, x_1, ..., x_{k-1} \in \PCTLsetAMB, x_k \in \PCTLsetB\}  \\ 
\end{aligned}
\end{equation}
where $\sqcup$ denotes a disjoint union.
The probability of the set
$\bigl\{ \PCTLBuntil{\PCTLsetA}{k}{\PCTLsetB} \bigr\}$
can be computed directly using the additivity of $\PP_x$:
\begin{equation}
\begin{aligned}
\label{UNTIL_DIRECT_COMPUTATION}
& \PP_{x} \left(\PCTLBuntil{\PCTLsetA}{k}{\PCTLsetB}\right)
= \\
&
\begin{cases}
1 & \text{if $x\in \PCTLsetB$} \\
\PP_{x}(x_1 \in \PCTLsetB) + \PP_x(x_1\in\PCTLsetAMB, x_2\in\PCTLsetB) +\cdots & \text{if $x\in \PCTLsetAMB$}\\
	\qquad \cdots + \PP_{x}(x_1, ..., x_{k-1} \in \PCTLsetAMB, x_k \in \PCTLsetB) & \\
0 & {\rm otherwise}
\end{cases}
\end{aligned}
\end{equation}
%
%
%
By the Markov property, all the latter probabilities can be expressed in terms of $Q$. For instance:
\begin{displaymath}
\begin{split}
& \PP_{x} \left( x_1, ..., x_{k-1} \in \PCTLsetAMB, x_k \in \PCTLsetB \right) \\
&\quad =
\int_{\PCTLsetAMB} Q(x,       \mrm d\xi_1) \cdots \int_{\PCTLsetAMB} Q(\xi_{k-3}, \mrm d\xi_{k-2}) \cdot\\
& \qquad\qquad\int_{\PCTLsetAMB} Q(\xi_{k-2}, \mrm d\xi_{k-1}) Q(\xi_{k-1}, \PCTLsetB).
\end{split}
\end{displaymath}

Nevertheless, $\PP_{x} \left(\PCTLBuntil{\PCTLsetA}{k}{\PCTLsetB}\right)$
can be computed more expressively in a recursive fashion.
Let $\MeasurableBounded(\StateSpaceX)$ be the set of all the measurable and bounded
functions defined over $\StateSpaceX$. $\MeasurableBounded(\StateSpaceX)$ is
a Banach space with the norm $\left\lVert f \right\rVert_\infty := \sup_{x\in\StateSpaceX}f(x)$.
Let the operator $\Lop:\MeasurableBounded(\StateSpaceX)\rightarrow \MeasurableBounded(\StateSpaceX)$ be defined as follows:
\begin{equation}
\label{e:Lop}
\Lop[W](x) := \Char{\PCTLsetB}(x) + \Char{\PCTLsetAMB}(x) \int_{\StateSpaceX} Q(x, \mrm d\xi) W(\xi)
\end{equation}
Given $\PCTLsetA$ and $\PCTLsetB$, 
let $\SolutionsSet \subset \MeasurableBounded(\StateSpaceX)$ be the set of functions
$W$ such that:
\begin{itemize}[label=$\circ$, leftmargin=*, itemsep=0pt, topsep=0pt, partopsep=0pt]
\item for all $x\in\StateSpaceX$, $0\leq W(x) \leq 1$;
\item for all $x\in\PCTLsetB$, $W(x)=1$;
\item for all $x\in \StateSpaceX \backslash (\PCTLsetA \cup \PCTLsetB)$, $W(x)=0$.
\end{itemize}

\begin{lemma}
\label{L_MAPS_LEMMA}
The set $\SolutionsSet$ is closed in $\MeasurableBounded(\StateSpaceX)$, and
$\Lop$ maps $\SolutionsSet$ into itself.
\end{lemma}
\begin{proof}
The closedness of $\SolutionsSet$ is trivial, because 
all of its three defining property are preserved even by {\em pointwise}
convergence.
Let $W\in\SolutionsSet$. 
The measurability of $\Lop[W]$ follows from the fact that if 
$Q$ is a stochastic kernel, and $W$
is a measurable bounded function, then the function $x \mapsto \int_\StateSpaceX Q(x, \mrm d\xi) W(\xi)$
is also measurable and bounded (see for instance \cite[Appendix C]{ref:hernandez-lerma1}).
The bounds $0 \leq \Lop[W](x) \leq 1$ are obvious,
since the same bounds hold for the integral, $Q(x, \cdot)$ being a probability on $\StateSpaceX$.
The fact that $\Lop[W](x)=1 \ \forall x\in\PCTLsetB$ and $\Lop[W](x) = 0 \ \forall x\in \StateSpaceX \backslash (\PCTLsetA \cup \PCTLsetB)$
is also obvious due to the indicator functions in the definition of $\Lop$.
\end{proof}

For fixed $\PCTLsetA$ and $\PCTLsetB$, 
let us now define recursively:
\begin{equation}
\label{VK_DEF}
\begin{split}
V_0     &:= \Char{\PCTLsetB} \\
V_{k+1} &:= \Lop[V_k]
\end{split}
\end{equation}

\begin{lemma}
\label{VK_EQUALTO_LEMMA}
\label{VK_INCREASING_LEMMA}
For all $k\geq 0$, 
$V_k(x) \equiv \PP_{x} \left( \PCTLBuntil{\PCTLsetA}{k}{\PCTLsetB} \right)$.
Moreover, for all $x$, the sequence $\{V_k(x)\}$ is nondecreasing.
\end{lemma}

\begin{proof}
Substituting recursively $V_1$ into $V_2$, $V_2$ into $V_3$ and so on, we obtain
\begin{align*}
V_2(x) &= 
          \Char{\PCTLsetB}(x) 
        + \Char{\PCTLsetAMB}(x) Q(x, \PCTLsetB)\\
		& \quad\;\; + \Char{\PCTLsetAMB}(x)\int_\PCTLsetAMB Q(x, d\xi_1) Q(\xi_1, \PCTLsetB),
        \\
V_3(x) &=
          \Char{\PCTLsetB}(x) 
        + \Char{\PCTLsetAMB}(x)     Q(x, \PCTLsetB)\\
		& \quad\;\; + \Char{\PCTLsetAMB}(x)     \int_\PCTLsetAMB Q(x, d\xi_1) Q(\xi_1, \PCTLsetB) \\
        & \quad\;\; + \Char{\PCTLsetAMB}(x)\int_\PCTLsetAMB Q(x, d\xi_1) \int_\PCTLsetAMB Q(\xi_1, d\xi_2) Q(\xi_2, \PCTLsetB)
        \\
	    &...\\
V_k(x) &=
          \Char{\PCTLsetB}(x) 
        + \Char{\PCTLsetAMB}(x) \cdot\\
		& \sum_{i=1}^{k} \overset{i-1 \  \text{times}}{\overbrace{\int_\PCTLsetAMB Q(x, d\xi_1) \cdots \int_\PCTLsetAMB Q(\xi_{i-2}, d\xi_{i-1}) }} Q(\xi_{i-1}, \PCTLsetB).
\end{align*}
Then, by the Markov property,
\begin{displaymath}
\begin{split}
V_k(x) &= \Char{\PCTLsetB}(x) 
        + \Char{\PCTLsetAMB}(x)\! \sum_{i=1}^{k} \PP_{x}( x_1, ..., x_{i-1} \in \PCTLsetAMB, x_i \in \PCTLsetB) \\
       &= \PP_{x} \left( x\in\PCTLsetB \right) 
        + \sum_{i=1}^{k} \PP_{x} \left( x, x_1, ..., x_{i-1} \in \PCTLsetAMB, x_i \in \PCTLsetB \right) \\
       &= \PP_{x} \left( \{x\in \PCTLsetB\}    \sqcup ... \sqcup \{x, x_1, ..., x_{k-1} \in \PCTLsetAMB, x_k \in \PCTLsetB\} \right) \\ 
       &= \PP_{x} \left( \PCTLBuntil{\PCTLsetA}{k}{\PCTLsetB} \right)\\
\end{split}
\end{displaymath}
The first assertion is proved. 
The second one is easily proved by induction. 
Obviously $V_1(x) - V_0(x) = \Char{\PCTLsetAMB}(x) Q(x, \PCTLsetB) \geq 0$. Suppose now that
$V_{k+1}(x) - V_k(x) \geq 0$. Then
$V_{k+2}(x) - V_{k+1}(x) = 
\Char{\PCTLsetAMB}(x) \int_\StateSpaceX Q(x, \mrm d\xi) \left(V_{k+1}(\xi) - V_k(\xi) \right) \geq 0$.
It follows by induction that for all $k\ge 0$ and all $x\in \StateSpaceX$ we have $V_{k+1}(x) \geq V_k(x)$.
\end{proof}

The semantics of the ``bounded until'' PCTL operator is now easy to explain.
In view of Lemma \ref{VK_EQUALTO_LEMMA}, given $\PCTLsetA$ and $\PCTLsetB$ we have:
\begin{displaymath}
x \vDash  \PCTLPP{\PCTLBuntil{\PCTLsetA}{k}{\PCTLsetB}}  \quad \Leftrightarrow \quad V_k(x)  \sim p
\end{displaymath}
Since $V_k$ is Borel measurable, any super- or sub-level set of the kind
$\PCTLPP{\PCTLBuntil{\PCTLsetA}{k}{\PCTLsetB}}$
is a Borel subset of $\StateSpaceX$.

\subsection{``Unbounded until''}
\label{SUBSECTION_UNBOUNDED_UNTIL}

Finally, we develop the ``unbounded until'' PCTL formula in detail.
Suppose, as before, that the process starts from $x_0 = x$. 
In the process's probability space we consider the event
\begin{equation}
\label{UNBOUNDED_UNTIL_EVENT}
\begin{split}
&\left\{ \PCTLuntil{\PCTLsetA}{\PCTLsetB} \right\}
= \left\{ \exists \tau \in \Nz: x, x_1, ..., x_{\tau-1} \in \PCTLsetA, x_\tau \in \PCTLsetB \right\}  \\ 
&\;\;=  \{x \in \PCTLsetB\} \cup ...   \cup   \{x,x_1, ..., x_{k-1} \in \PCTLsetA,   x_k \in \PCTLsetB\} \cup ...   \\ 
&\;\;=  \{x \in \PCTLsetB\} \sqcup ... \sqcup \{x,x_1, ..., x_{k-1} \in \PCTLsetAMB, x_k \in \PCTLsetB\} \sqcup ... 
\end{split}
\end{equation}
Its probability is as follows:
\begin{equation}
\label{UNBOUNDED_UNTIL_DIRECT_COMPUTATION}
\begin{aligned}
&\PP_{x} \left(\PCTLuntil{\PCTLsetA}{\PCTLsetB}\right) = \\
&\begin{cases}
1 & \text{if $x\in \PCTLsetB$,} \\
\sum_{k=1}^{+\infty} \PP_{x}(x_1, ..., x_{k-1} \in \PCTLsetAMB, x_k \in \PCTLsetB) & \text{if $x\in\PCTLsetAMB$,}\\
0 &  \text{otherwise.}
\end{cases}
\end{aligned}
\end{equation}
%
%
Notice, however, that the ``unbounded until'' event
is indeed the limit of the nondecreasing sequence of ``bounded until''
events we have considered above, i.e.,
\begin{displaymath}
\left\{ \PCTLuntil{\PCTLsetA}{\PCTLsetB} \right\}
=
\bigcup_{k=0}^{+\infty}
\left\{ \PCTLBuntil{\PCTLsetA}{k}{\PCTLsetB}\right\}
\end{displaymath}
Consequently, for all $x$ its probability can be obtained as the following limit:
\begin{displaymath}
\PP_{x} \left(\PCTLuntil{\PCTLsetA}{\PCTLsetB}\right)
= \lim_{k\rightarrow+\infty} \PP_{x} \left(\PCTLBuntil{\PCTLsetA}{k}{\PCTLsetB}\right)
= \lim_{k\rightarrow+\infty} V_k(x)
\end{displaymath}
(This limit is also a {\em supremum}, since the $V_k$ form a nondecreasing sequence.)
We define
\begin{equation}
\label{V_DEFINITION}
V(x) := \lim_{k\rightarrow+\infty} V_k(x)
\end{equation}

\begin{lemma}
\label{LEMMA_V_EQUATION}
The function $V$ defined in (\ref{V_DEFINITION}) 
belongs to $\SolutionsSet$ and
satisfies the following integral equation:
\begin{equation}
\label{V_EQUATION}
V(x) = \Char{\PCTLsetB}(x) + \Char{\PCTLsetAMB}(x) \int_\StateSpaceX Q(x, \mrm d\xi) V(\xi)
\end{equation}
(In other words, it is a fixed point for $\Lop$.)
\end{lemma}

\begin{proof}
The three properties required for the belonging to $\SolutionsSet$
are immediate, for they hold for all the $V_k$'s.
Consider again the recursive definition (\ref{VK_DEF}):
\begin{equation}
\label{VK_DEF2}
V_{k+1}(x) = \Char{\PCTLsetB}(x) + \Char{\PCTLsetAMB}(x) \int_{\StateSpaceX} Q(x, \mrm d\xi) V_k(\xi)
\end{equation}

From Lemmas
\ref{L_MAPS_LEMMA} and \ref{VK_INCREASING_LEMMA}, 
the $V_k$'s are Borel measurable and non-negative, and they form
a nondecreasing sequence. By definition of $V$, they converge pointwise to $V$.
Therefore, by the monotone convergence theorem 
(see for instance \cite[Theorem~1, p.~13]{ref:raoProbTheo})
for all $x$ we have
\begin{displaymath}
\lim_{k\rightarrow +\infty} \int_\StateSpaceX Q(x, \mrm d\xi) V_k(\xi)
= \int_\StateSpaceX Q(x, \mrm d\xi) V(\xi)
\end{displaymath}
Hence, letting $k\rightarrow+\infty$ in both sides of (\ref{VK_DEF2}), we
obtain (\ref{V_EQUATION}).
\end{proof}

The semantics of the ``unbounded until'' PCTL operator is now obvious.
For given $\PCTLsetA$ and $\PCTLsetB$, we have:
\begin{displaymath}
x \vDash  \PCTLPP{\PCTLuntil{\PCTLsetA}{\PCTLsetB}}
 \quad \Leftrightarrow \quad  V(x) \sim p
\end{displaymath}
Since $V$ is the limit of measurable functions, it is measurable itself,
hence its super- or sub-level sets
$\PCTLPP{\PCTLuntil{\PCTLsetA}{\PCTLsetB}}$
are again Borel subsets of $\StateSpaceX$.

\subsection{Notes on equation (\ref{V_EQUATION})}
\label{SUBSECTION_EQUATION}

First of all, note that the function $V$ defined in (\ref{V_DEFINITION})
is indeed a solution to equation (\ref{V_EQUATION}),
but it is by no means guaranteed to be its {\em unique} solution.
As a counterexample, let us consider the operator $\Eventually$
we have mentioned in the finite case. Let $\PCTLsetB$ be a formula (set).
The path formula $\Eventually\PCTLsetB$ (``eventually $\PCTLsetB$'')
is defined as $\PCTLuntil{\True}{\PCTLsetB}$.
Its probability $V(x) = \PP_x\left(\PCTLuntil{\True}{\PCTLsetB}\right)$
must therefore satisfy:
\begin{equation}
\label{NON_UNIQUENESS_EXAMPLE}
V(x) = \Char{\PCTLsetB}(x) + \Char{\PCTLsetB^C}(x) \int_\StateSpaceX Q(x, \mrm d\xi) V(\xi)
\end{equation}
Suppose that the set $\PCTLsetB^C$ is absorbing (that is, $Q(x, \PCTLsetB) = 0$ for all $x\in\PCTLsetB^C$).
Then, it is easy to see that both $V(x) \equiv \Char{\PCTLsetB}(x)$
and $V(x) \equiv 1$ are solutions of (\ref{NON_UNIQUENESS_EXAMPLE})
(the meaningful one being the former).
As another limit example, consider
the event $\Eventually\False$ (``eventually, false will hold true''!).
Its probability, both by immediate intuition and by calculation,
must be zero for all $x$. Nevertheless, 
any constant function $V$ is a solution to the corresponding equation:
\begin{displaymath}
\label{NON_UNIQUENESS_EXAMPLE2}
V(x) = \Char{\varnothing}(x) + \Char{\StateSpaceX}(x) \int_\StateSpaceX Q(x, \mrm d\xi) V(\xi)
 = \int_\StateSpaceX Q(x, \mrm d\xi) V(\xi)
\end{displaymath}
(irrespective of the structure of $Q$).

We can get around this issue with a characterization of $V$
among the solutions of (\ref{V_EQUATION}). We have the following
result:
\begin{lemma}
Let $\{W_\alpha\}$ be the family of all the {\em non-negative} 
solutions to \eqref{V_EQUATION}, i.e.,
\begin{displaymath}
W_\alpha(x) = \Char{\PCTLsetB}(x) + \Char{\PCTLsetAMB}(x) \int_\StateSpaceX Q(x, \mrm d\xi) W_\alpha(\xi)
\end{displaymath}
Then, for all $x$
\begin{displaymath}
V(x) 
= \inf_{\alpha} W_\alpha(x)
\equiv \min_{\alpha} W_\alpha(x)
\end{displaymath}
\end{lemma}

\begin{proof}
First, we show that, for any $V_k$ defined in (\ref{VK_DEF}) and for any non-negative solution $W$
to (\ref{V_EQUATION}), it holds $V_k(x) \leq W(x)$.
%
%
Define $V_{-1}(x) \equiv 0$ on $\StateSpaceX$. Then we have
$\Lop[V_{-1}] = V_0$.
Now, for all $x\in\StateSpaceX$, $W(x) - V_{-1}(x) = W(x) \geq 0$ by hypothesis.
Assume that $W(x) - V_k(x) \geq 0$ for all $x$. Then
\begin{align*}
W(x) - V_{k+1}(x)
&= \Lop[W](x) - \Lop[V_k](x) \\
&= \Char{\PCTLsetAMB}(x) \int_{\StateSpaceX} Q(x, \mrm d\xi) \left( W(\xi) - V_k(\xi)\right) \\
&\geq 0.
\end{align*}
It follows by induction that $V_k(x) \leq W(x)$ for all $x\in\StateSpaceX$ and for all $k\in\Nz$.

Since the above inequality holds for all of the $V_k$'s, it also holds
for their supremum $V$, that is, $V(x) \leq W(x)$ for any non-negative solution $W$
to (\ref{V_EQUATION}). The assertion follows since, by Lemma \ref{LEMMA_V_EQUATION}, $V$ is itself a solution
to (\ref{V_EQUATION}).
\end{proof}

\section{Uniqueness of $V$}
\label{SECTION_STOPPING_TIMES}

\newcommand{\Lrestr}{\bar{L}}
\newcommand{\lnorm}{\left\lVert}
\newcommand{\rnorm}{\right\rVert}

This section treats the issue of uniqueness of solutions to the integral equation~\eqref{V_EQUATION}. We approach the problem from two different directions, the first is functional analytic:
\begin{proposition}
\label{CONTRACTION_PROPOSITION}
Suppose that
\begin{displaymath}
\sup_{x\in\PCTLsetAMB}Q(x,\PCTLsetAMB) < 1.
\end{displaymath}
Then 
\begin{enumerate}[label={\em (\arabic*)}, align=right, leftmargin=*]
\item $\Lop$ is a contraction on $\SolutionsSet$;
\item equation \ref{V_EQUATION} has a unique solution $V$;
\item the elements $V_k$ defined in \ref{VK_DEF} converge to $V$ in the $\lnorm\cdot\rnorm_\infty$ norm,
that is {\em uniformly} in $\StateSpaceX$.
\end{enumerate}
\end{proposition}

\begin{proof}
Let $\alpha = \sup_{x\in\PCTLsetAMB}Q(x,\PCTLsetAMB)$.
Let $W_1, W_2 \in\SolutionsSet$.
For all $x\in \PCTLsetB \cup \PCTLsetA^C$ $\left\lvert \Lop[W_1](x) - \Lop[W_2](x) \right\rvert = 0$, whereas for all $x\in \PCTLsetAMB$, we have
\begin{align*}
& \left\lvert \Lop[W_1](x) - \Lop[W_2](x) \right\rvert\\
&\quad= \left\lvert \int_{\StateSpaceX} Q(x, \mrm d\xi) W_1(\xi) - \int_{\StateSpaceX} Q(x, \mrm d\xi) W_2(\xi) \right\rvert \\
&\quad\leq \int_{\StateSpaceX} Q(x, \mrm d\xi) \ \left\lvert W_1(\xi) - W_2(\xi) \right\rvert \\
&\quad= \int_{\PCTLsetAMB} Q(x, \mrm d\xi) \ \left\lvert W_1(\xi) - W_2(\xi) \right\rvert \\
&\quad\leq \int_{\PCTLsetAMB} Q(x, \mrm d\xi) \ \left\lVert W_1 - W_2 \right\rVert_\infty \\
&\quad= Q(x, \PCTLsetAMB) \ \left\lVert W_1 - W_2 \right\rVert_\infty \\ 
&\quad\leq \alpha \  \left\lVert W_1 - W_2 \right\rVert_\infty
\end{align*}
%
Since the above bound holds for each $x$, it holds also for the supremum over
$\PCTLsetAMB$, and consequently for the supremum over $\StateSpaceX$:
\begin{align*}
\left\lVert \Lop[W_1] - \Lop[W_2] \right\rVert_\infty 
& = \sup_{x\in\StateSpaceX} \left\lvert \Lop[W_1](x) - \Lop[W_2](x) \right\rvert \\
& 
\leq \alpha \  \left\lVert W_1 - W_2 \right\rVert_\infty
\end{align*}

This concludes the proof of claim (1). Claims (2) and (3) follow
by the Contraction Mapping Theorem 
\cite[Theorem 9.23]{ref:rudinPrinciples}
since $\SolutionsSet$ is closed.
\end{proof}

\begin{corollary}
Suppose that $\sup_{x\in\PCTLsetAMB}Q(x,\PCTLsetAMB) < 1$.
Suppose moreover that $Q$ satisfies the strong Feller or strong continuity~\cite[Appendix C]{ref:hernandez-lerma1} property.
Then the restriction of $V$ to $\PCTLsetAMB$ is continuous.
\end{corollary}


\begin{proof}
Let $\bar{V}$ and $\bar{V}_k$ denote the restriction to $\PCTLsetAMB$
of $V$ and $V_k$ respectively. In particular, we have
\begin{equation}
\label{VK_RESTRICTED}
\begin{split}
\bar{V}_0(x) &= 0 \\
\bar{V}_{k+1}(x) 
&= \int_\StateSpaceX Q(x, \mrm d\xi) V_k(\xi)\\
& \qquad = \int_\PCTLsetAMB Q(x, \mrm d\xi) \bar{V}_k(\xi) + Q(x, \PCTLsetB) 
\end{split}
\end{equation}
Obviously $\bar{V}_0$ is continuous.  Due to the strong Feller property, $x \mapsto Q(x, \PCTLsetB)$ is continuous, and if $\bar{V}_k$ is measurable then $x \mapsto \int_\PCTLsetAMB Q(x, \mrm d\xi) \bar{V}_k(\xi)$ and therefore $\bar{V}_{k+1}$ are continuous.  By induction, all the $\bar{V}_k$ are continuous.  Hence, $\{\bar{V}_k\}$ is a sequence of continuous functions that converges uniformly to $\bar{V}$.  Thus, $\bar{V}$ is also continuous.
\end{proof}

The second direction is probabilistic: Let us define two random times 
	\begin{equation}
	\label{e:taudef}
	\begin{aligned}
		\tau & := \inf\bigl\{t\in\Nz\,\big|\, X_t\in \PCTLsetB\bigr\} \;\;\text{and}\\
		\tau' & := \inf\bigl\{t\in\Nz\,\big|\, X_t\in X\setminus(\Phi\cup\Psi)\bigr\}.
	\end{aligned}
	\end{equation}
	It is not difficult to see that $\tau$ and $\tau'$ are stopping times with respect to the filtration $(\sigalg_t)_{t\in\Nz}$. Also observe that
	\begin{align*}
		\PP_x\bigl(\Phi\Until\Psi\bigr) & = \PP_x(\tau < \tau', \tau < \infty),\\
		V(x) & = \PP_x(\tau < \tau', \tau < \infty) = \EE_x\!\Biggl[\sum_{t=0}^{\tau\mn\tau'}\Char{\Psi}(X_t)\Biggr].
	\end{align*}

\begin{proposition}
	Assume that $\tau\mn\tau' < \infty$ almost surely. Then, for $u\in\SolutionsSet$ we have
	\begin{enumerate}[label={\em (\roman*)}, leftmargin=*, widest=iii]
		\item $u\le V$ whenever $u$ satisfies the functional inequality $u\le \Lop[u]$, and
		\item $u\ge V$ whenever $u$ satisfies $u\ge \Lop[u]$,
	\end{enumerate}
	where all inequalities are interpreted pointwise on $X$. In particular, $V$ is the unique solution to the equation $u = \Lop[u]$ on the set $\SolutionsSet$.
\end{proposition}
\begin{proof}
	We prove (i) first.	Fix $u\in\SolutionsSet$ and $x\in X$. From Lemma \ref{VK_INCREASING_LEMMA} it follows readily that $\Lop$ is a monotone operator on $\MeasurableBounded$. Iterating the inequality $u\le \Lop[u]$ $n$-times we arrive at
	\begin{align*}
		& u(x) \le \Lop[u](x)\\
			& \le \Lop[\Lop[u]](x) \le\cdots \le \overset{n-\text{times}}{\overbrace{\Lop[\Lop[\cdots\Lop}}[u]\cdots]]\\
			& = \Char{\Psi}(x) + \Char{\Phi\setminus\Psi}(x)\int_{\Phi\cup\Psi}Q(x, \mrm d\xi_1)\biggl(\indic{\Psi}(\xi_1) + \\
			& \qquad \indic{\Phi\setminus\Psi}(\xi_1)\int_{\Phi\cup\Psi}Q(\xi_1, \mrm d\xi_2)\biggl(\cdots + \cdots \biggl(\indic{\Psi}(\xi_{n-1}) + \\
			& \qquad \indic{\Phi\setmin\Psi}(\xi_{n-1})\int_{\Phi\cup\Psi}Q(\xi_{n-1}, \mrm d\xi_n)u(\xi_n)\biggr)\biggr)\biggr)\\
			& = \Biggl(\Char{\Psi}(x) + \Char{\Phi\setminus\Psi}(x)\int_{\Phi\cup\Psi}Q(x, \mrm d\xi_1)\biggl(\indic{\Psi}(\xi_1) + \\
			& \quad \indic{\Phi\setminus\Psi}(\xi_1)\int_{\Phi\cup\Psi}Q(\xi_1, \mrm d\xi_2)\biggl(\cdots + \cdots \biggl(\indic{\Psi}(\xi_{n-2}) + \\
			& \quad \indic{\Phi\setmin\Psi}(\xi_{n-2})\int_{\Phi\cup\Psi}Q(\xi_{n-2}, \mrm d\xi_{n-1})\indic{\Psi}(\xi_{n-1})\biggr)\cdots\biggr)\biggr)\Biggr)\\
			& \quad + \Biggl(\indic{\Phi\setmin\Psi}(x)\int_{\Phi\setmin\Psi}Q(x, \mrm d\xi_1)\cdots\int_{\Phi\setmin\Psi}Q(\xi_{n-2}, \mrm d\xi_{n-1})\cdot \\
			& \qquad\qquad\qquad\qquad\int_{\Phi\cup\Psi}Q(\xi_{n-1}, \mrm d\xi_{n})u(\xi_n)\Biggr)\\
			& = \EE_x\Biggl[\sum_{t=0}^{(n-1)\mn\tau\mn\tau'}\indic{\Psi}(X_t)\Biggr] + \\
			& \;\;\quad \EE_x\!\left[\indic{\Phi\setmin\Psi}(X_{(n-1)\mn\tau\mn\tau'})(\indic{\Phi\cup\Psi}\cdot u)(X_{n\mn\tau\mn\tau'})\indic{\{\tau\mn\tau' < \infty\}}\right].
	\end{align*}
	The left-hand side above is independent of $n$, and since $\tau\mn\tau' < \infty$ almost surely, taking limits we get
	\begin{align*}
		& u(x) \le \lim_{n\to\infty}\EE_x\Biggl[\sum_{t=0}^{(n-1)\mn\tau\mn\tau'}\indic{\Psi}(X_t)\Biggr] +\\
		& \qquad \lim_{n\to\infty}\EE_x\bigl[\indic{\Phi\setmin\Psi}(X_{(n-1)\mn\tau\mn\tau'})\cdot\\
		& \qquad\qquad\qquad\qquad(\indic{\Phi\cup\Psi}\cdot u)(X_{n\mn\tau\mn\tau'})\indic{\{\tau\mn\tau' < \infty\}}\bigr]\\
		& \quad = \EE_x\Biggl[\sum_{t=0}^{\tau\mn\tau'}\indic{\Psi}(X_t)\Biggr] +\\
		& \qquad \EE_x\bigl[\indic{\Phi\setmin\Psi}(X_{\tau\mn\tau'})(\indic{\Phi\cup\Psi}\cdot u)(X_{\tau\mn\tau'})\indic{\{\tau\mn\tau' < \infty\}}\bigr]\\
		& \quad = V(x) + 0.
	\end{align*}
	To justify the interchange of integration and limit above we have employed the monotone and the dominated convergence theorems for the first and the second terms, respectively, and since $X_{\tau\mn\tau'}\not\in \Phi\setmin\Psi$ by definition, the last expectation vanishes. Since $u\in\SolutionsSet$ and $x\in X$ are arbitrary, we see that $u\le\Lop[u]$ implies $u\le V$ whenever $u\in\SolutionsSet$. The proof of (ii) follows exactly the same arguments as above, with ``$\ge$'' replacing every ``$\le$'' everywhere in the above steps; we omit the details. Uniqueness of $V$ as a solution of the functional equation $u = \Lop[u]$ on the set $\SolutionsSet$ follows at once from (i) and (ii).
\end{proof}

\section{Examples}
\label{SECTION_EXAMPLES}

We demonstrate the effectiveness of the PCTL verification methodology on two simple problems with potentially important implications.  The first example comes from the literature on fishery management, where multiple recovery strategies for a single species fishery are considered.  The second example comes from the finance literature, where the problem of early retirement is explored.  In both examples, the problems are solved numerically by gridding the state space.  It is of great interest to pursue more effective and accurate solution methods for the DP integral equations using sophisticate methods commented on in the Introduction.

\subsection{Recovery Strategies in Fishery Management}
Overexploitation can lead to both a decrease in the fish stock to a level below which maximum sustainable yield (MSY) cannot be supported and/or a decrease in fish stock to a level where net revenue has been driven to zero \cite{Clark}.  When the fish stock drops below this level, appropriate recovery strategies are necessary to recover the fish stock while minimizing economic loss.  In this example, we use the PCTL framework to evaluate the effectiveness of various recovery strategies (or non-strategies) over a finite time horizon for the recovery of a fish population.

We consider a discrete time Markov model of a single species fishery motivated by \cite{Pitchford}.  For a time horizon $k = 0,1,\ldots N$, the evolution of the fish biomass within a fishable area is given by the stochastic difference equation \cite{Pitchford}
\[
x_{k+1} = (1 - \nu_k)x_k + \gamma_kR(x_k) - \delta_kC(x_k),
\]
where $x_k$ is the fish biomass at time $k$, $R(\cdot)$ is a function representing the recruitment (e.g., addition through birth) of fish, $C(\cdot)$ is the catch function, $\nu_k$ is a random variable that represents fish mortality during stage $k$, $\gamma_k$ is a random variable representing the variability in the recruitment of the fish population, and $\delta_k$ is a random variable representing the variability in the catch.  The species recruitment function is given by
\[
R(x_k) = \max\biggl\{rx_k\left(1 - \frac{x_k}{2K}\right),0\biggr\},
\]
where $r\in [0,1]$ is the per-capita recruitment at time step $k$ and $K$ is equal to half the biomass limit (i.e., upper bound on the fish population) for the fishable area. 

We consider three different recovery strategies implemented through the target catch function.  In the first, we apply a constant target catch according to the deterministic MSY \cite{Kot}, i.e.
\[
C(x_k) = C_{\text{MSY}} = \frac{K(r - \mu)^2}{2r},
\]
where $\mu$ is the deterministic mortality rate.  The second recovery strategy is given by the Harvest Control Rule (HCR)
\[
C(x_k) = \begin{cases}
	C_{\text{MSY}}\frac{x_k}{K} & \text{if } x_k < K,\\
	C_{\text{MSY}} & \text{otherwise.} 
\end{cases}
\]
Lastly, we consider the strategy $C(x_k) = 0$.

Following \cite{Pitchford}, we assign the values $K = 200$, $r = 1$, and $\mu = 0.2$, and take all random variables to be i.i.d. according to the following distributions $\nu \sim \mathcal{N}(\mu,0.1^2)$,  $\gamma \sim \mathcal{N}(1,0.6^2)$, and $\delta \sim \mathcal{N}(1.1,0.2^2)$.  Using the MSY as a measure of safety for the system, we assign the target operating region for the fishery to be $K = [150, 400]$ and the safe operating region to be $K' = \;]0,400]$.

\begin{figure*}[ht]
\centering
\subfigure[Value Function]{\includegraphics[keepaspectratio=true,width=5cm]{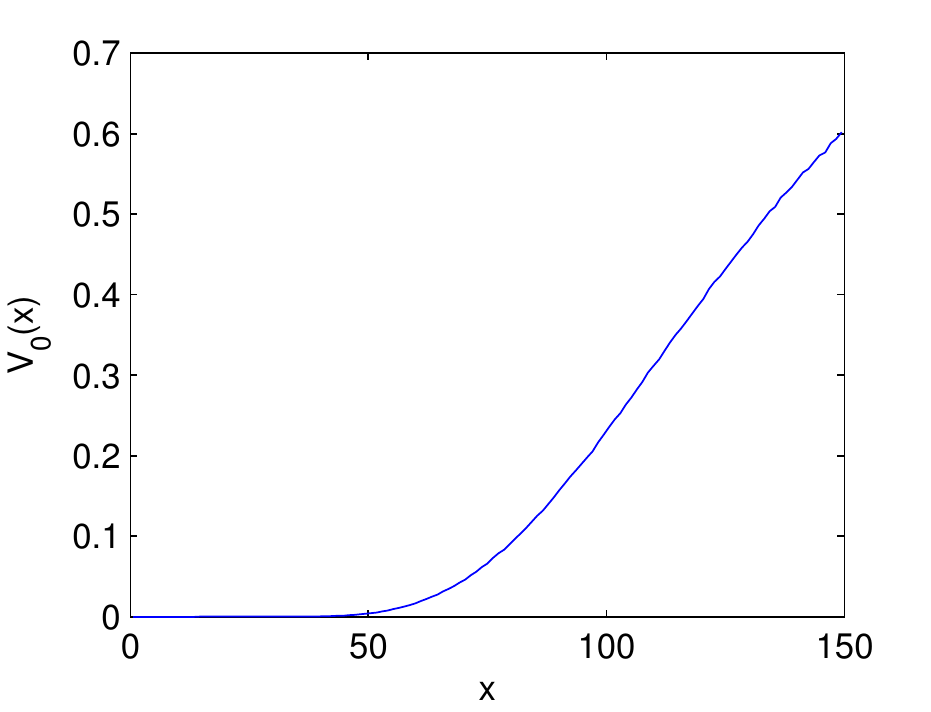}}
\subfigure[Value Function]{\includegraphics[keepaspectratio=true,width=5cm]{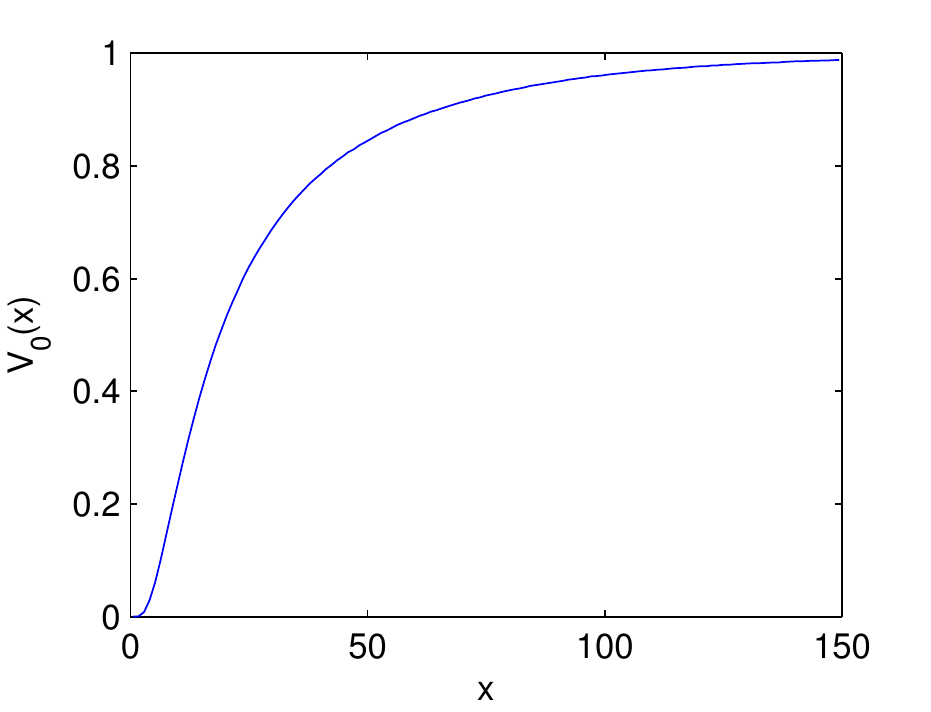}}
\subfigure[Value Function]{\includegraphics[keepaspectratio=true,width=5cm]{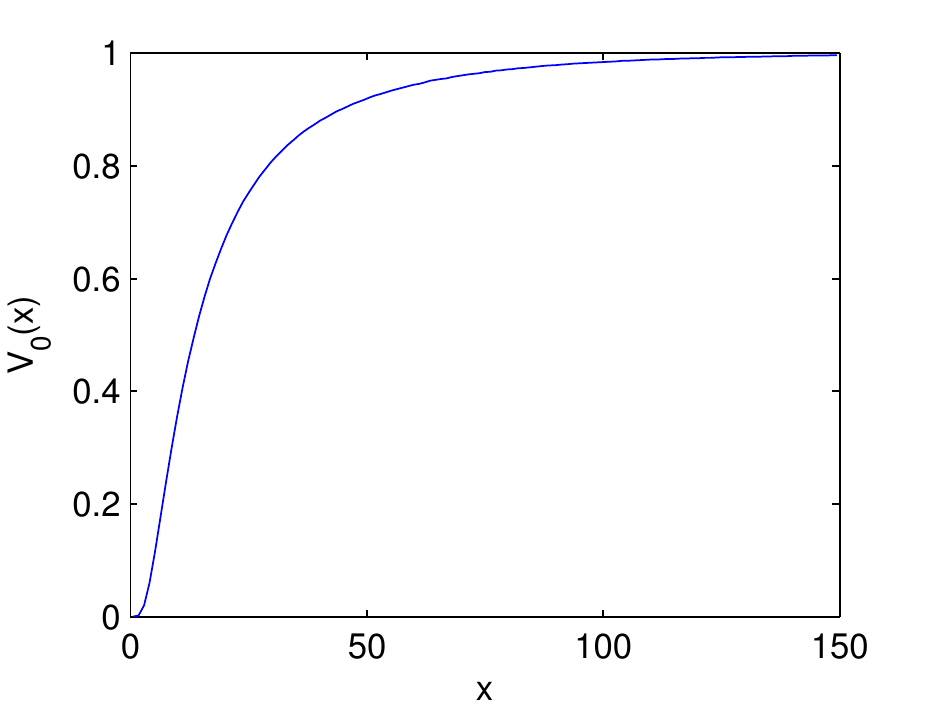}}
\caption{Results for the Recovery Problem at time $k = 0$.  The function $V_0(\cdot)$ for different recovery policies are given in (a) MSY, (b) HCR, and (c) Fishing Stop.
}\label{Recover}
\end{figure*}

For the verification of the control strategies, we consider the set of initial states (i.e., fish biomass at $k = 0$) that satisfy
\[
P_{\geq 0.9}\bigl[K' \mathcal{U}^{\leq 5} K\bigr].
\]
That is, we are interested in the set of states that, with a probability greater than $90$ percent, will enter the target operating region $K$ within $N = 5$ time steps while remaining in $K'$ until then.  The functions satisfying the dynamic recursion \eqref{VK_DEF2} for the three different recovery strategies are shown in Figure \ref{Recover}.  According to the computational results, the sets that satisfy the bounded until operator are approximately $\emptyset$, $[65,400]$, and $[45,400]$ for the three policies respectively.  It is interesting to note that under the deterministic MSY quota policy the solution is the empty set, meaning that there are no initial states which result in recovery with $90$ percent certainty over the short time horizon.  Further, the gain in reliable recovery between the HCR strategy and a complete fishing stop is minimal, indicating that it may be in the economic interest of the fishery to use the HCR policy in the region.

%
\subsection{A Problem of Early Retirement}


Recently, increased attention has been given to stochastic risk models with investment income in the discrete time setting \cite{Cai1,deKok,Wei,Nyr,Tang1,Cai2,Yang}.  In most cases the probability of ruin over a finite or infinite time horizon is the main area of interest, with the infinite horizon case being mathematically easier and thus more popular in the literature \cite{Paulson2}.  Interestingly enough, personal retirement funds fall into the same category as basic ruin models, and therefore can be modeled as such.  Further, the individual is often as concerned with the short term financial gain (e.g., achieving a financial target for the fund) as with the risk of losing the investment (i.e., ruin).

Motivated by \cite{Target}, we consider a discrete time Markov model of an individual retirement fund.  Based on \cite{Paulson2}, the evolution of the retirement fund $x_k$ over a finite horizon $k = 0,1,\ldots,N$ is given according to the stochastic difference equation
\[
x_{k+1} = a x_k(1 + S_k) + b x_k(1 + R_k) + c x_k + u_k,
\]
where $x_k$ is the value of the retirement fund and $u_k$ is the yearly individual contribution to the fund.  $S_k$ and $R_k$ are i.i.d. random variables representing the average rates of return for a safe investment and a risky investment over one year, $a$ is the percentage of capital invested in the safe asset, $b$ is the percentage of capital invested in the risky asset, and $c$ is the percentage capital not invested at all.  Note the restriction that $a + b + c = 1$.

For simplicity, all random variables are assumed to be i.i.d. with $S_k \sim \mathcal{N}(0.03,0.005^2)$ and $R_k \sim \mathcal{N}(0.1,0.2^2)$ for all $k = 0,1,\ldots,N$.  We consider three different investment strategies (i) $a = 0.4$, $b = 0.4$, and $c = 0.2$, (ii) $a = 0.8$, $b = 0.2$, and $c = 0$, and (iii) $a = 0.2$, $b = 0.8$, and $c = 0$.  For each strategy, the yearly contribution is $u_k = 2500$ for all $k = 0,1,\ldots,N$.

Consider the target set $K = [200000, +\infty[$ and the safe set $K' = \;]0,+\infty[$.  Over a finite time horizon of $N = 20$ years, we would like to identify the set of all initial investments $x_0 \in \R$ such that the retirement fund hits the target set $K$ (i.e., surpasses $200000$) while avoiding total financial ruin with a probability greater than $85$ percent.  To this end, we consider the PCTL formula
\begin{eqnarray}
P_{\geq 0.85}\bigl[K' \mathcal{U}^{\leq 20} K\bigr].
\end{eqnarray}

\begin{figure*}[ht]
\centering
\subfigure[Value Function]{\includegraphics[keepaspectratio=true,width=5cm]{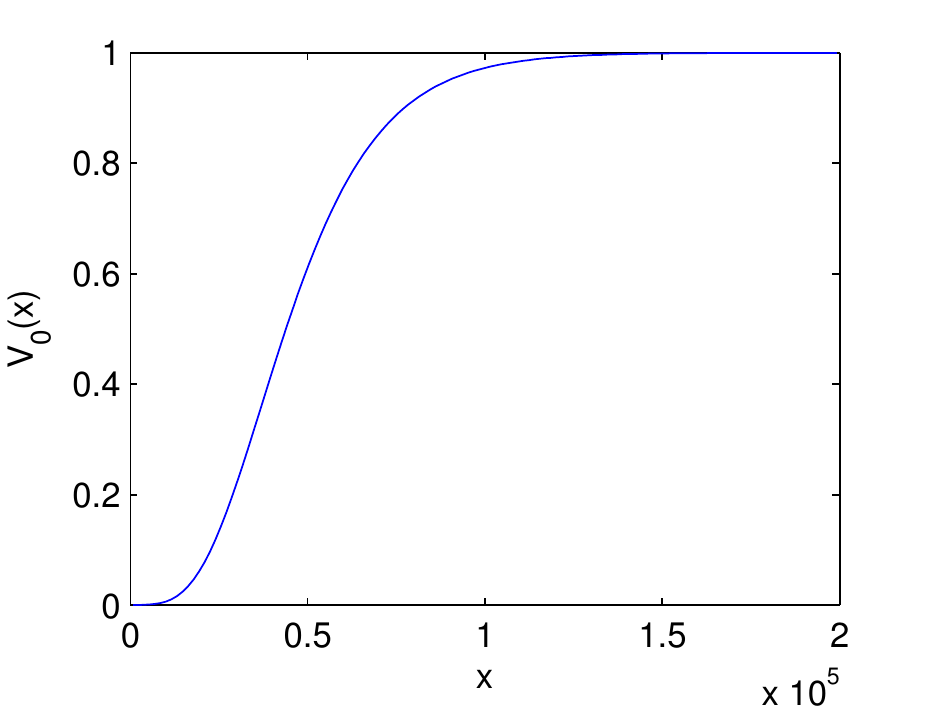}}
\subfigure[Value Function]{\includegraphics[keepaspectratio=true,width=5cm]{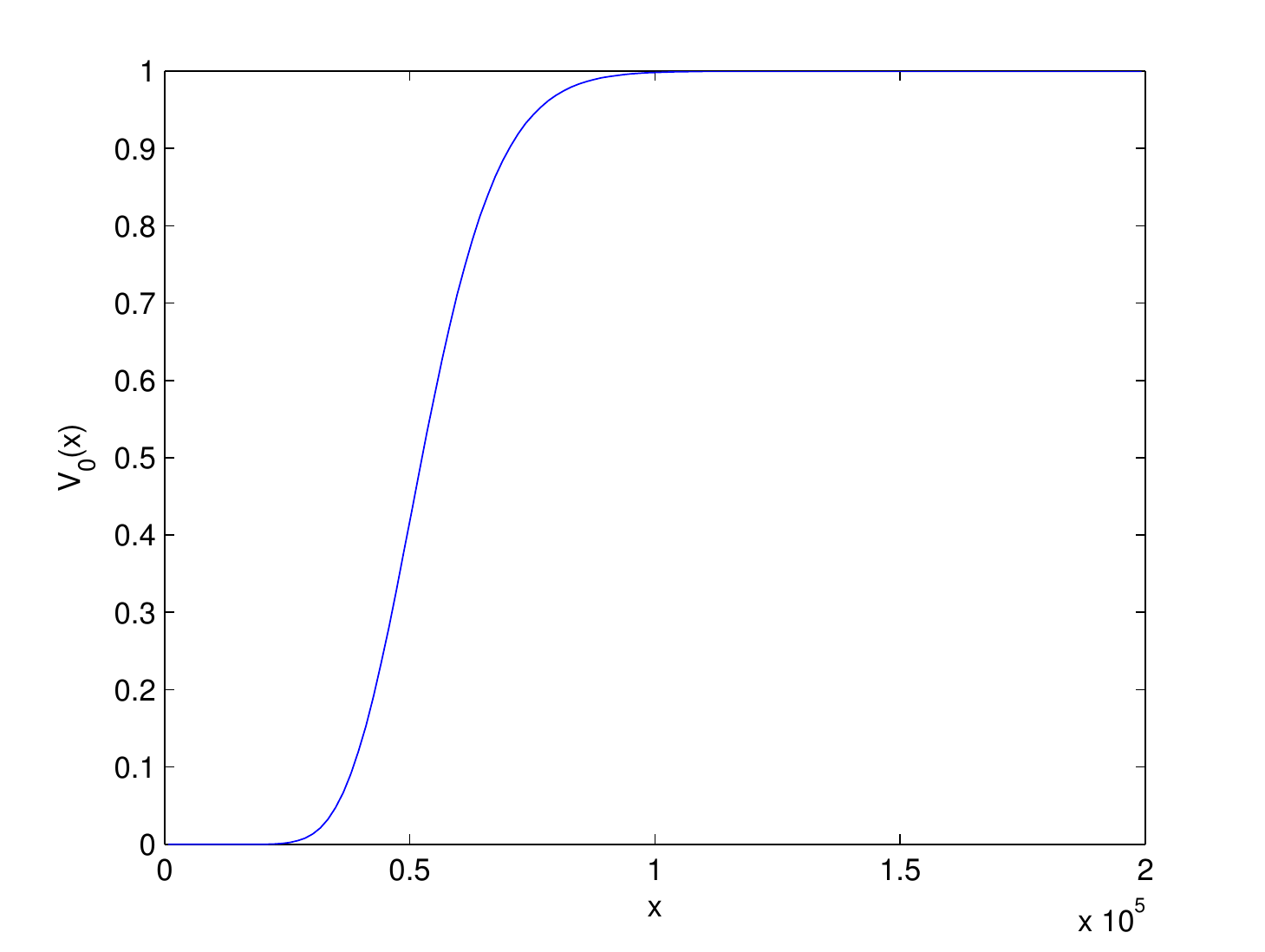}}
\subfigure[Value Function]{\includegraphics[keepaspectratio=true,width=5cm]{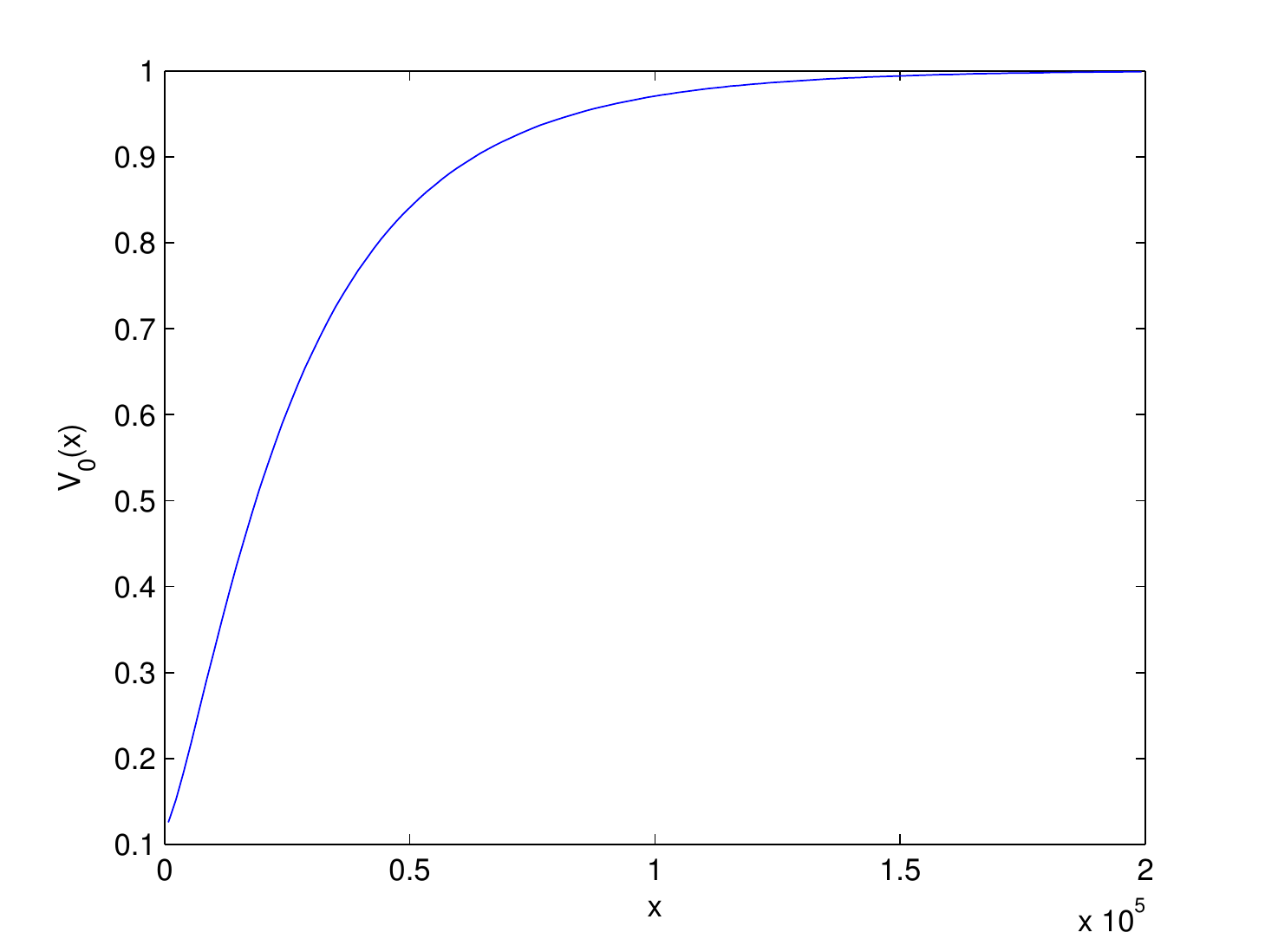}}
\caption{The function $V_0(\cdot)$ for the Early Retirement problem at year $k = 0$ for different investment policies are given in (a) $a=0.4,\enspace b=0.4,\enspace c=0.2$, (b) $a=0.8,\enspace b=0.2,\enspace c=0$, and (c) $a=0.2,\enspace b=0.8, \enspace c=0$.
}\label{Retire}
\end{figure*}

For each investment strategy, the function satisfying the dynamic recursion \eqref{VK_DEF2} at time $k=0$ is shown in Figure \ref{Retire}.  According to the computational results, the set that satisfies the bounded until operator for each strategy is given by (i) $[70000, +\infty[$, (ii) $[66500, +\infty[$, and (iii) $[51500, +\infty[$.  Thus, with an initial investment of more than $51500$ swiss francs, yearly contributions in the amount of $2500$ swiss francs, and investment strategy (iii), an individual has an $85$ percent chance of retiring within $20$ years.  However, if we were to consider an increasing probability of success, at some point the strategy with the largest set satisfying the bounded until operator would switch from (iii) to (ii).

\section{Conclusion and Future Work}

In this paper, we have extended the grammar and semantics of PCTL for {\em finite}-state Markov chains for the verification of general state-space Markov chains.  We have shown that the bulk of the computational methodology is in the evaluation of the ``bounded until'' and ``unbounded until'' operators.  And that the evaluation of these operators reduces to the computation of DP-like integral equations, for which there is a rich numerical history.

In the future, extensions to the language to capture additional trajectories will be explored which maintain the DP-like structure.  Also, numerical methods for the efficient and accurate evaluation of the DP integral equations are being evaluated and applied to various sample problems.


\end{document}